\documentclass{gtart}

\def\ifplaintex{\expandafter\ifx\csname documentclass\endcsname\relax}

\def\gtp{{\mathsurround=0pt\it $\cal G\mskip-2mu$eometry \&\ 
$\cal T\!\!$opology $\cal P\!$ublications}}  

\def\recd{{\small Received:\qua\receiveddate\ifx\reviseddate\relax
\else\qquad Revised:\qua\reviseddate\fi\par}} 

\def\lognumber#1{\def\thelognumber{#1}}
\def\volumenumber#1{\def\thevolumenumber{#1}}
\def\volumeyear#1{\def\thevolumeyear{#1}}
\def\papernumber#1{\def\thepapernumber{#1}}
\def\pagenumbers#1#2{\def\startpage{#1}\def\finishpage{#2}}
\def\published#1{\def\publishdate{#1}}

\def\received#1{\def\receiveddate{#1}}
\def\revised#1{\def\reviseddate{#1}}
\def\accepted#1{\def\accepteddate{#1}}

\def\asciiauthors#1{\def\theasciiauthors{#1}}
\def\asciiaddress#1{\def\theasciiaddress{#1}}

\def\coverauthors#1{\def\thecoverauthors{#1}}

\let\\\par\let\thelognumber\relax\let\thevolumenumber\relax
\let\thepapernumber\relax\let\thevolumeyear\relax\let\startpage\relax
\let\finishpage\relax\let\publishdate\relax\let\receiveddate\relax
\let\reviseddate\relax\let\accepteddate\relax\let\theasciititle\relax
\let\theasciiauthors\relax\let\theasciiaddress\relax
\let\theasciiabstract\relax

\let\thecoverauthors\relax\let\theasciiemail\relax

\ifplaintex
\font\logobig=cmssbx10 scaled 3836
\font\logomed=cmssbx10 scaled 2557
\else
\font\logobig=cmssbx10 scaled 4200
\font\logomed=cmssbx10 scaled 2800
\fi

\long\def\makeagttitle{   
\count0=\startpage
\agt\hfill      
\hbox to 45truept{\vbox to 0pt{\vglue -13truept{\logomed A\kern -.37em{\logobig 
T}\kern -.38em G}\vss}\hss}
\break
{\small Volume \thevolumenumber\ (\thevolumeyear)
\startpage--\finishpage\nl
Published: \publishdate}

\vglue .25truein

{\parskip=0pt\leftskip 0pt plus
1fil\def\\{\par\smallskip}{\Large\bf\thetitle}\par\medskip} \vglue
0.05truein

{\parskip=0pt\leftskip 0pt plus 1fil\def\\{\par}{\sc\theauthors}
\par\medskip}%
 
\vglue 0.03truein 

{\small\leftskip 25truept\rightskip 25truept{\bf Abstract}\stdspace\theabstract

{\bf AMS Classification}\stdspace\theprimaryclass
\ifx\thesecondaryclass\relax\else; \thesecondaryclass\fi\par
{\bf Keywords}\stdspace \thekeywords\par}\vglue 7truept

}   

\ifplaintex
\font\phead=cmsl9 scaled 950
\font\pnum=cmbx10 scaled 913
\font\pfoot=cmsl9 scaled 950
\headline{\vbox to 0pt{\vskip -4.5mm\line{\small\phead\ifnum
\count0=\startpage ISSN 1472-2739 (on-line) 1472-2747 (printed)
\hfill {\pnum\folio}\else\ifodd\count0\def\\{ }%
\ifx\theshorttitle\relax\thetitle\else\theshorttitle\fi\hfill{\pnum\folio}
\else\def\\{ and }{\pnum\folio}\hfill\ifx\theshortauthors\relax\theauthors
\else\theshortauthors\fi\fi\fi}\vss}}
\footline{\vbox to 0pt{\vglue 0mm\line{\small\pfoot\ifnum\count0=\startpage
\copyright\ \gtp\hfill\else
\agt, Volume \thevolumenumber\ (\thevolumeyear)\hfill\fi}\vss}}
\else
\font=cmsl9 scaled 1050
\font\lnum=cmbx10 
\font=cmsl9 scaled 1050
\makeatletter
\def\@oddhead{{\small\lhead\ifnum\count0=\startpage ISSN 1472-2739 
(on-line) 1472-2747 (printed)\hfill {\lnum\number\count0}\else\ifodd\count0
\def\\{ }\ifx\theshorttitle\relax \thetitle \else\theshorttitle\fi\hfill
{\lnum\number\count0}\else\def\\{ and }{\lnum\number\count0}
\hfill\ifx\theshortauthors\relax 
\theauthors\else\theshortauthors\fi\fi\fi}}\def\@evenhead{\@oddhead}
\def\@oddfoot{\small\lfoot\ifnum\count0=\startpage\copyright\ \gtp\hfill\else
\agt, Volume \thevolumenumber\ (\thevolumeyear)\hfill\fi}
\def\@evenfoot{\@oddfoot}
\makeatother
\fi
\let\maketitlepage\makeagttitle
\let\makeshorttitle\maketitlepage
\let\maketitle\maketitlepage

\newwrite\gtoutfile
\long\gdef\makeheadfile{  
{\def\\{, }\def\s{ }
\immediate\openout\gtoutfile head.xxx
\immediate\write\gtoutfile{To: math@arxiv.org}
\immediate\write\gtoutfile{Subject: put OR rep NNNNN:ppppp}
\immediate\write\gtoutfile{--text follows this line--}
\immediate\write\gtoutfile{Proxy-for: \ifx\theasciiauthors\relax
\theauthors\else\theasciiauthors\fi\s<\ifx\theasciiemail\relax\theemail\else\theasciiemail\fi>}
\immediate\write\gtoutfile{\noexpand\\}
\immediate\write\gtoutfile{Authors: \ifx\theasciiauthors\relax
\theauthors\else\theasciiauthors\fi}
{\def\\{ }\immediate\write\gtoutfile{Title: \ifx\theasciititle\relax
\thetitle\else\theasciititle\fi}}
\immediate\write\gtoutfile{Subj-class: GT or SG, GR etc}
\immediate\write\gtoutfile{MSC-class: \theprimaryclass\ifx\thesecondaryclass\relax\else, \thesecondaryclass\fi}
\immediate\write\gtoutfile{Journal-ref: Algebr. Geom. Topol. \thevolumenumber\s
(\thevolumeyear) \startpage-\finishpage}
\immediate\write\gtoutfile{Comments: Published by Algebraic and
Geometric Topology at}
\immediate\write\gtoutfile{\s\s\s  http://www.maths.warwick.ac.uk/agt/AGTVol\thevolumenumber/agt-\thevolumenumber-\thepapernumber.abs.html}
\immediate\write\gtoutfile{\noexpand\\}
\immediate\write\gtoutfile{}
\ifx\theasciiabstract\relax
\immediate\write\gtoutfile{\theabstract}\else
\immediate\write\gtoutfile{\theasciiabstract}\fi
\immediate\write\gtoutfile{}
\immediate\write\gtoutfile{\noexpand\\}
\immediate\write\gtoutfile{}
\immediate\closeout\gtoutfile}}  

\def\maketitlepage{\makeagttitle\makeheadfile}
\let\makeshorttitle\maketitlepage
\let\maketitle\maketitlepage

\def\ifplaintex{\expandafter\ifx\csname documentclass\endcsname\relax}

\def\gtp{{\mathsurround=0pt\it $\cal G\mskip-2mu$eometry \&\ 
$\cal T\!\!$opology $\cal P\!$ublications}}  

\def\recd{{\small Received:\qua\receiveddate\ifx\reviseddate\relax
\else\qquad Revised:\qua\reviseddate\fi\par}} 

\def\lognumber#1{\def\thelognumber{#1}}
\def\volumenumber#1{\def\thevolumenumber{#1}}
\def\volumeyear#1{\def\thevolumeyear{#1}}
\def\papernumber#1{\def\thepapernumber{#1}}
\def\pagenumbers#1#2{\def\startpage{#1}\def\finishpage{#2}}
\def\published#1{\def\publishdate{#1}}

\def\received#1{\def\receiveddate{#1}}
\def\revised#1{\def\reviseddate{#1}}
\def\accepted#1{\def\accepteddate{#1}}

\def\asciiauthors#1{\def\theasciiauthors{#1}}
\def\asciiaddress#1{\def\theasciiaddress{#1}}

\def\coverauthors#1{\def\thecoverauthors{#1}}

\let\\\par\let\thelognumber\relax\let\thevolumenumber\relax
\let\thepapernumber\relax\let\thevolumeyear\relax\let\startpage\relax
\let\finishpage\relax\let\publishdate\relax\let\receiveddate\relax
\let\reviseddate\relax\let\accepteddate\relax\let\theasciititle\relax
\let\theasciiauthors\relax\let\theasciiaddress\relax
\let\theasciiabstract\relax

\let\thecoverauthors\relax\let\theasciiemail\relax

\ifplaintex
\font\logobig=cmssbx10 scaled 3836
\font\logomed=cmssbx10 scaled 2557
\else
\font\logobig=cmssbx10 scaled 4200
\font\logomed=cmssbx10 scaled 2800
\fi

\long\def\makeagttitle{   
\count0=\startpage
\agt\hfill      
\hbox to 45truept{\vbox to 0pt{\vglue -13truept{\logomed A\kern -.37em{\logobig 
T}\kern -.38em G}\vss}\hss}
\break
{\small Volume \thevolumenumber\ (\thevolumeyear)
\startpage--\finishpage\nl
Published: \publishdate}

\vglue .25truein

{\parskip=0pt\leftskip 0pt plus
1fil\def\\{\par\smallskip}{\Large\bf\thetitle}\par\medskip} \vglue
0.05truein

{\parskip=0pt\leftskip 0pt plus 1fil\def\\{\par}{\sc\theauthors}
\par\medskip}%
 
\vglue 0.03truein 

{\small\leftskip 25truept\rightskip 25truept{\bf Abstract}\stdspace\theabstract

{\bf AMS Classification}\stdspace\theprimaryclass
\ifx\thesecondaryclass\relax\else; \thesecondaryclass\fi\par
{\bf Keywords}\stdspace \thekeywords\par}\vglue 7truept

}   

\ifplaintex
\font\phead=cmsl9 scaled 950
\font\pnum=cmbx10 scaled 913
\font\pfoot=cmsl9 scaled 950
\headline{\vbox to 0pt{\vskip -4.5mm\line{\small\phead\ifnum
\count0=\startpage ISSN 1472-2739 (on-line) 1472-2747 (printed)
\hfill {\pnum\folio}\else\ifodd\count0\def\\{ }%
\ifx\theshorttitle\relax\thetitle\else\theshorttitle\fi\hfill{\pnum\folio}
\else\def\\{ and }{\pnum\folio}\hfill\ifx\theshortauthors\relax\theauthors
\else\theshortauthors\fi\fi\fi}\vss}}
\footline{\vbox to 0pt{\vglue 0mm\line{\small\pfoot\ifnum\count0=\startpage
\copyright\ \gtp\hfill\else
\agt, Volume \thevolumenumber\ (\thevolumeyear)\hfill\fi}\vss}}
\else
\font=cmsl9 scaled 1050
\font\lnum=cmbx10 
\font=cmsl9 scaled 1050
\makeatletter
\def\@oddhead{{\small\lhead\ifnum\count0=\startpage ISSN 1472-2739 
(on-line) 1472-2747 (printed)\hfill {\lnum\number\count0}\else\ifodd\count0
\def\\{ }\ifx\theshorttitle\relax \thetitle \else\theshorttitle\fi\hfill
{\lnum\number\count0}\else\def\\{ and }{\lnum\number\count0}
\hfill\ifx\theshortauthors\relax 
\theauthors\else\theshortauthors\fi\fi\fi}}\def\@evenhead{\@oddhead}
\def\@oddfoot{\small\lfoot\ifnum\count0=\startpage\copyright\ \gtp\hfill\else
\agt, Volume \thevolumenumber\ (\thevolumeyear)\hfill\fi}
\def\@evenfoot{\@oddfoot}
\makeatother
\fi
\let\maketitlepage\makeagttitle
\let\makeshorttitle\maketitlepage
\let\maketitle\maketitlepage

\newwrite\gtoutfile
\long\gdef\makeheadfile{  
{\def\\{, }\def\s{ }
\immediate\openout\gtoutfile head.xxx
\immediate\write\gtoutfile{To: math@arxiv.org}
\immediate\write\gtoutfile{Subject: put OR rep NNNNN:ppppp}
\immediate\write\gtoutfile{--text follows this line--}
\immediate\write\gtoutfile{Proxy-for: \ifx\theasciiauthors\relax
\theauthors\else\theasciiauthors\fi\s<\ifx\theasciiemail\relax\theemail\else\theasciiemail\fi>}
\immediate\write\gtoutfile{\noexpand\\}
\immediate\write\gtoutfile{Authors: \ifx\theasciiauthors\relax
\theauthors\else\theasciiauthors\fi}
{\def\\{ }\immediate\write\gtoutfile{Title: \ifx\theasciititle\relax
\thetitle\else\theasciititle\fi}}
\immediate\write\gtoutfile{Subj-class: GT or SG, GR etc}
\immediate\write\gtoutfile{MSC-class: \theprimaryclass\ifx\thesecondaryclass\relax\else, \thesecondaryclass\fi}
\immediate\write\gtoutfile{Journal-ref: Algebr. Geom. Topol. \thevolumenumber\s
(\thevolumeyear) \startpage-\finishpage}
\immediate\write\gtoutfile{Comments: Published by Algebraic and
Geometric Topology at}
\immediate\write\gtoutfile{\s\s\s  http://www.maths.warwick.ac.uk/agt/AGTVol\thevolumenumber/agt-\thevolumenumber-\thepapernumber.abs.html}
\immediate\write\gtoutfile{\noexpand\\}
\immediate\write\gtoutfile{}
\ifx\theasciiabstract\relax
\immediate\write\gtoutfile{\theabstract}\else
\immediate\write\gtoutfile{\theasciiabstract}\fi
\immediate\write\gtoutfile{}
\immediate\write\gtoutfile{\noexpand\\}
\immediate\write\gtoutfile{}
\immediate\closeout\gtoutfile}}  

\def\maketitlepage{\makeagttitle\makeheadfile}
\let\makeshorttitle\maketitlepage
\let\maketitle\maketitlepage

\def\ifplaintex{\expandafter\ifx\csname documentclass\endcsname\relax}

\def\gtp{{\mathsurround=0pt\it $\cal G\mskip-2mu$eometry \&\ 
$\cal T\!\!$opology $\cal P\!$ublications}}  

\def\recd{{\small Received:\qua\receiveddate\ifx\reviseddate\relax
\else\qquad Revised:\qua\reviseddate\fi\par}} 

\def\lognumber#1{\def\thelognumber{#1}}
\def\volumenumber#1{\def\thevolumenumber{#1}}
\def\volumeyear#1{\def\thevolumeyear{#1}}
\def\papernumber#1{\def\thepapernumber{#1}}
\def\pagenumbers#1#2{\def\startpage{#1}\def\finishpage{#2}}
\def\published#1{\def\publishdate{#1}}

\def\received#1{\def\receiveddate{#1}}
\def\revised#1{\def\reviseddate{#1}}
\def\accepted#1{\def\accepteddate{#1}}

\def\asciiauthors#1{\def\theasciiauthors{#1}}
\def\asciiaddress#1{\def\theasciiaddress{#1}}

\def\coverauthors#1{\def\thecoverauthors{#1}}

\let\\\par\let\thelognumber\relax\let\thevolumenumber\relax
\let\thepapernumber\relax\let\thevolumeyear\relax\let\startpage\relax
\let\finishpage\relax\let\publishdate\relax\let\receiveddate\relax
\let\reviseddate\relax\let\accepteddate\relax\let\theasciititle\relax
\let\theasciiauthors\relax\let\theasciiaddress\relax
\let\theasciiabstract\relax

\let\thecoverauthors\relax\let\theasciiemail\relax

\ifplaintex
\font\logobig=cmssbx10 scaled 3836
\font\logomed=cmssbx10 scaled 2557
\else
\font\logobig=cmssbx10 scaled 4200
\font\logomed=cmssbx10 scaled 2800
\fi

\long\def\makeagttitle{   
\count0=\startpage
\agt\hfill      
\hbox to 45truept{\vbox to 0pt{\vglue -13truept{\logomed A\kern -.37em{\logobig 
T}\kern -.38em G}\vss}\hss}
\break
{\small Volume \thevolumenumber\ (\thevolumeyear)
\startpage--\finishpage\nl
Published: \publishdate}

\vglue .25truein

{\parskip=0pt\leftskip 0pt plus
1fil\def\\{\par\smallskip}{\Large\bf\thetitle}\par\medskip} \vglue
0.05truein

{\parskip=0pt\leftskip 0pt plus 1fil\def\\{\par}{\sc\theauthors}
\par\medskip}%
 
\vglue 0.03truein 

{\small\leftskip 25truept\rightskip 25truept{\bf Abstract}\stdspace\theabstract

{\bf AMS Classification}\stdspace\theprimaryclass
\ifx\thesecondaryclass\relax\else; \thesecondaryclass\fi\par
{\bf Keywords}\stdspace \thekeywords\par}\vglue 7truept

}   

\ifplaintex
\font\phead=cmsl9 scaled 950
\font\pnum=cmbx10 scaled 913
\font\pfoot=cmsl9 scaled 950
\headline{\vbox to 0pt{\vskip -4.5mm\line{\small\phead\ifnum
\count0=\startpage ISSN 1472-2739 (on-line) 1472-2747 (printed)
\hfill {\pnum\folio}\else\ifodd\count0\def\\{ }%
\ifx\theshorttitle\relax\thetitle\else\theshorttitle\fi\hfill{\pnum\folio}
\else\def\\{ and }{\pnum\folio}\hfill\ifx\theshortauthors\relax\theauthors
\else\theshortauthors\fi\fi\fi}\vss}}
\footline{\vbox to 0pt{\vglue 0mm\line{\small\pfoot\ifnum\count0=\startpage
\copyright\ \gtp\hfill\else
\agt, Volume \thevolumenumber\ (\thevolumeyear)\hfill\fi}\vss}}
\else
\font=cmsl9 scaled 1050
\font\lnum=cmbx10 
\font=cmsl9 scaled 1050
\makeatletter
\def\@oddhead{{\small\lhead\ifnum\count0=\startpage ISSN 1472-2739 
(on-line) 1472-2747 (printed)\hfill {\lnum\number\count0}\else\ifodd\count0
\def\\{ }\ifx\theshorttitle\relax \thetitle \else\theshorttitle\fi\hfill
{\lnum\number\count0}\else\def\\{ and }{\lnum\number\count0}
\hfill\ifx\theshortauthors\relax 
\theauthors\else\theshortauthors\fi\fi\fi}}\def\@evenhead{\@oddhead}
\def\@oddfoot{\small\lfoot\ifnum\count0=\startpage\copyright\ \gtp\hfill\else
\agt, Volume \thevolumenumber\ (\thevolumeyear)\hfill\fi}
\def\@evenfoot{\@oddfoot}
\makeatother
\fi
\let\maketitlepage\makeagttitle
\let\makeshorttitle\maketitlepage
\let\maketitle\maketitlepage

\newwrite\gtoutfile
\long\gdef\makeheadfile{  
{\def\\{, }\def\s{ }
\immediate\openout\gtoutfile head.xxx
\immediate\write\gtoutfile{To: math@arxiv.org}
\immediate\write\gtoutfile{Subject: put OR rep NNNNN:ppppp}
\immediate\write\gtoutfile{--text follows this line--}
\immediate\write\gtoutfile{Proxy-for: \ifx\theasciiauthors\relax
\theauthors\else\theasciiauthors\fi\s<\ifx\theasciiemail\relax\theemail\else\theasciiemail\fi>}
\immediate\write\gtoutfile{\noexpand\\}
\immediate\write\gtoutfile{Authors: \ifx\theasciiauthors\relax
\theauthors\else\theasciiauthors\fi}
{\def\\{ }\immediate\write\gtoutfile{Title: \ifx\theasciititle\relax
\thetitle\else\theasciititle\fi}}
\immediate\write\gtoutfile{Subj-class: GT or SG, GR etc}
\immediate\write\gtoutfile{MSC-class: \theprimaryclass\ifx\thesecondaryclass\relax\else, \thesecondaryclass\fi}
\immediate\write\gtoutfile{Journal-ref: Algebr. Geom. Topol. \thevolumenumber\s
(\thevolumeyear) \startpage-\finishpage}
\immediate\write\gtoutfile{Comments: Published by Algebraic and
Geometric Topology at}
\immediate\write\gtoutfile{\s\s\s  http://www.maths.warwick.ac.uk/agt/AGTVol\thevolumenumber/agt-\thevolumenumber-\thepapernumber.abs.html}
\immediate\write\gtoutfile{\noexpand\\}
\immediate\write\gtoutfile{}
\ifx\theasciiabstract\relax
\immediate\write\gtoutfile{\theabstract}\else
\immediate\write\gtoutfile{\theasciiabstract}\fi
\immediate\write\gtoutfile{}
\immediate\write\gtoutfile{\noexpand\\}
\immediate\write\gtoutfile{}
\immediate\closeout\gtoutfile}}  

\def\maketitlepage{\makeagttitle\makeheadfile}
\let\makeshorttitle\maketitlepage
\let\maketitle\maketitlepage

\def\ifplaintex{\expandafter\ifx\csname documentclass\endcsname\relax}

\def\gtp{{\mathsurround=0pt\it $\cal G\mskip-2mu$eometry \&\ 
$\cal T\!\!$opology $\cal P\!$ublications}}  

\def\recd{{\small Received:\qua\receiveddate\ifx\reviseddate\relax
\else\qquad Revised:\qua\reviseddate\fi\par}} 

\def\lognumber#1{\def\thelognumber{#1}}
\def\volumenumber#1{\def\thevolumenumber{#1}}
\def\volumeyear#1{\def\thevolumeyear{#1}}
\def\papernumber#1{\def\thepapernumber{#1}}
\def\pagenumbers#1#2{\def\startpage{#1}\def\finishpage{#2}}
\def\published#1{\def\publishdate{#1}}

\def\received#1{\def\receiveddate{#1}}
\def\revised#1{\def\reviseddate{#1}}
\def\accepted#1{\def\accepteddate{#1}}

\def\asciiauthors#1{\def\theasciiauthors{#1}}
\def\asciiaddress#1{\def\theasciiaddress{#1}}

\def\coverauthors#1{\def\thecoverauthors{#1}}

\let\\\par\let\thelognumber\relax\let\thevolumenumber\relax
\let\thepapernumber\relax\let\thevolumeyear\relax\let\startpage\relax
\let\finishpage\relax\let\publishdate\relax\let\receiveddate\relax
\let\reviseddate\relax\let\accepteddate\relax\let\theasciititle\relax
\let\theasciiauthors\relax\let\theasciiaddress\relax
\let\theasciiabstract\relax

\let\thecoverauthors\relax\let\theasciiemail\relax

\ifplaintex
\font\logobig=cmssbx10 scaled 3836
\font\logomed=cmssbx10 scaled 2557
\else
\font\logobig=cmssbx10 scaled 4200
\font\logomed=cmssbx10 scaled 2800
\fi

\long\def\makeagttitle{   
\count0=\startpage
\agt\hfill      
\hbox to 45truept{\vbox to 0pt{\vglue -13truept{\logomed A\kern -.37em{\logobig 
T}\kern -.38em G}\vss}\hss}
\break
{\small Volume \thevolumenumber\ (\thevolumeyear)
\startpage--\finishpage\nl
Published: \publishdate}

\vglue .25truein

{\parskip=0pt\leftskip 0pt plus
1fil\def\\{\par\smallskip}{\Large\bf\thetitle}\par\medskip} \vglue
0.05truein

{\parskip=0pt\leftskip 0pt plus 1fil\def\\{\par}{\sc\theauthors}
\par\medskip}%
 
\vglue 0.03truein 

{\small\leftskip 25truept\rightskip 25truept{\bf Abstract}\stdspace\theabstract

{\bf AMS Classification}\stdspace\theprimaryclass
\ifx\thesecondaryclass\relax\else; \thesecondaryclass\fi\par
{\bf Keywords}\stdspace \thekeywords\par}\vglue 7truept

}   

\ifplaintex
\font\phead=cmsl9 scaled 950
\font\pnum=cmbx10 scaled 913
\font\pfoot=cmsl9 scaled 950
\headline{\vbox to 0pt{\vskip -4.5mm\line{\small\phead\ifnum
\count0=\startpage ISSN 1472-2739 (on-line) 1472-2747 (printed)
\hfill {\pnum\folio}\else\ifodd\count0\def\\{ }%
\ifx\theshorttitle\relax\thetitle\else\theshorttitle\fi\hfill{\pnum\folio}
\else\def\\{ and }{\pnum\folio}\hfill\ifx\theshortauthors\relax\theauthors
\else\theshortauthors\fi\fi\fi}\vss}}
\footline{\vbox to 0pt{\vglue 0mm\line{\small\pfoot\ifnum\count0=\startpage
\copyright\ \gtp\hfill\else
\agt, Volume \thevolumenumber\ (\thevolumeyear)\hfill\fi}\vss}}
\else
\font=cmsl9 scaled 1050
\font\lnum=cmbx10 
\font=cmsl9 scaled 1050
\makeatletter
\def\@oddhead{{\small\lhead\ifnum\count0=\startpage ISSN 1472-2739 
(on-line) 1472-2747 (printed)\hfill {\lnum\number\count0}\else\ifodd\count0
\def\\{ }\ifx\theshorttitle\relax \thetitle \else\theshorttitle\fi\hfill
{\lnum\number\count0}\else\def\\{ and }{\lnum\number\count0}
\hfill\ifx\theshortauthors\relax 
\theauthors\else\theshortauthors\fi\fi\fi}}\def\@evenhead{\@oddhead}
\def\@oddfoot{\small\lfoot\ifnum\count0=\startpage\copyright\ \gtp\hfill\else
\agt, Volume \thevolumenumber\ (\thevolumeyear)\hfill\fi}
\def\@evenfoot{\@oddfoot}
\makeatother
\fi
\let\maketitlepage\makeagttitle
\let\makeshorttitle\maketitlepage
\let\maketitle\maketitlepage

\newwrite\gtoutfile
\long\gdef\makeheadfile{  
{\def\\{, }\def\s{ }
\immediate\openout\gtoutfile head.xxx
\immediate\write\gtoutfile{To: math@arxiv.org}
\immediate\write\gtoutfile{Subject: put OR rep NNNNN:ppppp}
\immediate\write\gtoutfile{--text follows this line--}
\immediate\write\gtoutfile{Proxy-for: \ifx\theasciiauthors\relax
\theauthors\else\theasciiauthors\fi\s<\ifx\theasciiemail\relax\theemail\else\theasciiemail\fi>}
\immediate\write\gtoutfile{\noexpand\\}
\immediate\write\gtoutfile{Authors: \ifx\theasciiauthors\relax
\theauthors\else\theasciiauthors\fi}
{\def\\{ }\immediate\write\gtoutfile{Title: \ifx\theasciititle\relax
\thetitle\else\theasciititle\fi}}
\immediate\write\gtoutfile{Subj-class: GT or SG, GR etc}
\immediate\write\gtoutfile{MSC-class: \theprimaryclass\ifx\thesecondaryclass\relax\else, \thesecondaryclass\fi}
\immediate\write\gtoutfile{Journal-ref: Algebr. Geom. Topol. \thevolumenumber\s
(\thevolumeyear) \startpage-\finishpage}
\immediate\write\gtoutfile{Comments: Published by Algebraic and
Geometric Topology at}
\immediate\write\gtoutfile{\s\s\s  http://www.maths.warwick.ac.uk/agt/AGTVol\thevolumenumber/agt-\thevolumenumber-\thepapernumber.abs.html}
\immediate\write\gtoutfile{\noexpand\\}
\immediate\write\gtoutfile{}
\ifx\theasciiabstract\relax
\immediate\write\gtoutfile{\theabstract}\else
\immediate\write\gtoutfile{\theasciiabstract}\fi
\immediate\write\gtoutfile{}
\immediate\write\gtoutfile{\noexpand\\}
\immediate\write\gtoutfile{}
\immediate\closeout\gtoutfile}}  

\def\maketitlepage{\makeagttitle\makeheadfile}
\let\makeshorttitle\maketitlepage
\let\maketitle\maketitlepage

\lognumber{34}
\volumenumber{2}
\volumeyear{2002}
\papernumber{34}
\published{8 October 2002}
\pagenumbers{825}{842}
\received{19 November 2001}
\revised{3 September 2002}
\accepted{6 September 2002}

\usepackage{amsmath,amssymb,epsfig}

\newtheorem{theorem}{Theorem}
\newtheorem{lemma}[theorem]{Lemma}
\newtheorem{corollary}[theorem]{Corollary}

\theoremstyle{definition}

\newtheorem{definition}[theorem]{Definition} 

\newtheorem{remark}[theorem]{Remark} 

\newcommand{\bR}{{\hspace{-0.3pt}\mathbb R}\hspace{0.3pt}}
\newcommand{\bZ}{{\hspace{-0.3pt}\mathbb Z}\hspace{0.3pt}}
\newcommand{\bC}{{\hspace{-0.3pt}\mathbb C}\hspace{0.3pt}}

\newcommand{\lcr}{\raisebox{-5pt}{\mbox{}\hspace{1pt} \
\epsfig{file=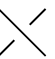} \hspace{1pt}\mbox{}}}
\newcommand{\ift}{\raisebox{-5pt}{\mbox{}\hspace{1pt}
 \ \epsfig{file=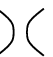} \hspace{1pt}\mbox{}}}
\newcommand{\zer}{\raisebox{-5pt}{\mbox{}\hspace{1pt}
 \ \epsfig{file=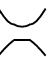} \hspace{1pt}\mbox{}}} 
\newcommand{\ga}{\gamma}
\newcommand{\de}{\delta} 
\newcommand{\ra}{\rightarrow} 

\newcommand{\<}{\langle}
\renewcommand{\>}{\rangle}

\begin{document}
\title{3-manifold invariants and periodicity\\of homology spheres}

\authors{Patrick M. Gilmer\\Joanna Kania-Bartoszynska\\J\'{o}zef H. Przytycki}
\shortauthors{Gilmer, Kania-Bartoszynska and Przytycki}

\coverauthors{Patrick M. Gilmer\\Joanna Kania-Bartoszynska\\J\noexpand\'{o}zef H. Przytycki}
\asciiauthors{Patrick M. Gilmer\\Joanna Kania-Bartoszynska\\Jozef H. Przytycki}
\email{gilmer@math.lsu.edu, kania@math.boisestate.edu, przytyck@research.circ.gwu.edu} 

\address{Department of Mathematics, Louisiana State University\\Baton 
Rouge, LA 70803, USA\\\smallskip\\Department of Mathematics, Boise State 
University\\Boise, ID 83725, USA\\\smallskip\\Department
of Mathematics, George Washington University\\Washington, D.C.20052, USA}

\asciiaddress{Department of Mathematics, Louisiana State University\\Baton 
Rouge, LA 70803, USA\\Department of Mathematics, Boise State 
University\\Boise, ID 83725, USA\\Department
of Mathematics, George Washington University\\Washington, D.C.20052, USA}



\begin{abstract} We show how the
periodicity of a homology sphere is reflected in the
Reshetikhin-Turaev-Witten invariants of the manifold. These yield a
criterion for the periodicity of a homology sphere. 
\end{abstract}

\maketitle 

{\small
\leftskip=25pt\rightskip=25pt\it

``\dots the 
angel of topology and the devil of abstract
  algebra fight for the soul of each individual mathematical domain" 
\rm\nl
\mbox{}\hfill Herman Weyl \cite{W}

\leftskip=0pt\rightskip=0pt}

\section{Introduction} 
Since 1988, we  have witnessed a
fascinating development in the theory of $3$-manifold invariants: from
Witten's ideas to Ohtsuki's finite type invariants. In this paper we
present one of the few nontrivial applications of this theory by showing
how the periodicity of a $3$-manifold is reflected in its quantum
invariants. 

In section \ref{2} we recall some necessary conditions for
periodicity of links in terms of the Jones polynomial and the Kauffman
bracket. Section \ref{3} contains Theorem \ref{periodicsurgery} 
which states that periodic homology spheres can be
presented as results of surgeries on periodic links. We also relate
periodic homology spheres to free cyclic covers of homology $S^1 \times
S^2$'s. In section \ref{4} we recall the definition and some properties of the
$SO(3)$-quantum invariants associated to a odd prime $p.$ These are 
versions of the 3-manifold invariants first described by Witten
\cite{Wi}, and shown to be well-defined by Reshetikhin-Turaev \cite{RT}.
Section \ref{5} contains  Theorem \ref{main}
which describes how the existence of a $\bZ_p$ action on a
homology sphere is reflected in the Reshetikhin-Turaev-Witten
invariants associated to  $p$-th roots of unity. 
Related results for $\bZ_p$ actions when the order of the root of unity
is co-prime to $p$ are  discussed in 
\cite{Gi1} and \cite{C}.  

In the last section we apply our criterion to the Brieskorn homology
spheres,  
$\Sigma(2 ,|n|, |2n-1|).$  In particular, we show that
the Poincar\'e homology sphere is not $p$-periodic where $p$ is any
prime between $7$ and $61$ (where we stopped the computer).

For the convenience of the reader, we include an  appendix written by
the first author with a direct proof of Theorem \ref{simplecong} . The
original proof  in \cite{Gi1} was given in a more general context.

Research on this paper was partially
supported by NSF-DMS-9626818, NSF-DMS-9971905, 
NSF-DMS-0203486 and by USAF grant
1-443964-22502.

\section{Periodic Links}\label{2} 

In this section we recall criteria for periodicity of links in $S^3$
in terms of the Jones polynomial and the Kauffman bracket. These criteria
motivated our main result, Theorem \ref{main}. 

Let $n$ be an integer greater than $1$.
\begin{definition}
  A link $L$ in $S^3$ is called $n$-periodic if there is a smooth
faithful $\bZ_n$-action on $S^3$, with a circle as 
the fixed point set, that maps $L$
onto itself, and where $L$ is disjoint from the fixed point set.
Furthermore, if $L$ is an oriented link, one assumes that a generator of
$\bZ_n$ preserves the orientation of $L$ or changes it to the opposite one.
\end{definition} 

\begin{figure}[ht!] 
\cl{\epsfig{file=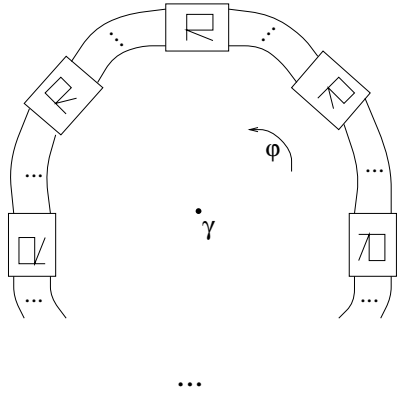,width=3in}\hspace{1in}\mbox{}}
\caption{Periodic link diagram} \label{perlink} 
\end{figure}

Given $n$-periodic link,
we know by the positive solution of the Smith Conjecture
\cite{MB}, 
that the fixed point set $\gamma$ of the $\bZ_n$
action  is an unknotted circle and this action is conjugate to an
orthogonal action on $S^3$. In other words, if we identify $S^3$ with
$\bR^3\cup\infty$, then the fixed point set can be assumed to be equal to
the ``vertical'' axis together with $\infty$. Then the generator of
$\bZ_n$ is the rotation $\varphi(z,t) = (e^{{2\pi i}/{n}}\cdot z,t)$, where
the coordinates on $\bR^3$ come from the product of the complex plane and
the real line, $\bC\times \bR$. The $n$-periodic link may be represented
by a $\varphi$-invariant diagram as in Figure \ref{perlink}.

\begin{figure}[ht!]\cl{\epsfig{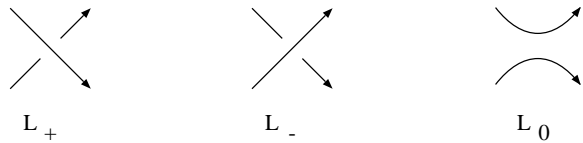}}
\caption{Skein triple} \label{Lplusminzer}
\end{figure}

Let $L_+, L_-, L_0$ denote three oriented links that are identical except in some 3-ball
as shown in Figure \ref{Lplusminzer}. The Jones polynomial \cite{J} of
oriented links, $V_L(t)\in Z[t^{\pm \frac{1}{2}}]$, is defined recursively
as follows: 

\begin{itemize} 
\item [(i)] $V_o = 1$, 
\item [(ii)]
$\frac{1}{t}V_{L_+}(t)- tV_{L_-}(t) = (\sqrt t - \frac{1}{\sqrt
  t})V_{L_0}(t)$. 
\end{itemize} 
The Kauffman bracket \cite{Ka} is a function on the set of unoriented
framed links in $\bR^3$. It takes values in the ring $\bZ[A^{\pm1}]$ and
is uniquely determined by the rules: 

\begin{enumerate} 
\item
$\displaystyle{\langle \underbrace{\bigcirc \dots \bigcirc}_m \rangle =
(-A^2-A^{-2})^{m-1}}$ \ \ and \vspace{5pt} 
\item
$\displaystyle{\left\langle\lcr\right\rangle =
                   A \left\langle\zer\right\rangle +
              A^{-1} \left\langle\ift\right\rangle}$. 
\end{enumerate} 
We assume, as usual, that the diagrams in condition (2) are identical in the
parts not shown, and that the links are given the blackboard framing. The
Kauffman bracket yields a simple definition of the Jones polynomial for
oriented links.  Let $L$ be an oriented link in $S^3$. For  $L$ equipped
 with some framing 
the writhe $w( L)$   (or {\em Tait number}) of $ L$ is the
linking number of the oriented link $L$ and the push-off of $L$ given by
the framing of $ L.$ 
 We have 
$$V_L(t) = (-A^3)^{-w( L ) }\langle L \rangle, $$ 
where $A = t^{-\frac{1}{4}}$.

One of the notable applications of polynomial knot
invariants is a criterion for periodicity discovered by Murasugi
\cite{M},  Przytycki \cite{P1} and Traczyk \cite{T}, as given in
the following theorem. 

\begin{theorem}\label{periodiclinks} 
Let $p$ be a prime number, and let $L$ be a
$p$-periodic oriented link in $S^3.$ Then $$V_L(t)\equiv
V_L(t^{-1}) \mod (p, t^p-1)$$ where $V_L(t)$ denotes Jones polynomial in
variable $t$, and $(p,t^p-1)$ is the ideal in $\bZ[t^{\pm \frac{1}{2}}]$
generated by $p$ and $t^p-1$.
\end{theorem} 

One obtains a  similar criterion involving the Kauffman bracket.

\begin{corollary}\label{periodicbracket} 
Let $ L$ be a $p$-periodic oriented  link where $p$ is prime. Assume $L$ has a framing which need not be periodic. We have that $$\< L\>_A \equiv A^{6 w( L) }
\<L\>_{A^{-1}} \mod{(p, A^{4p}-1)},$$ 
where we indicate by a
subscript that we change the variable from $A$ to $A^{-1}.$
\end{corollary}

Recall that a framed link is an embedding of a disjoint collection of
annuli into $S^3$. In diagrams we assume blackboard framing. We say
that a framed link is $n$-periodic if it has an $n$-periodic diagram.
In particular, an $n$-periodic framed knot must have writhe divisible
by $n$.

Since $w(L)\equiv 0 \pmod{p}$ for a $p$-periodic oriented framed link,  
we have 
\begin{corollary}\label{periodicbracket2} 
If $p$ is prime, $L$ is a $p$-periodic  
(unoriented)  framed link, and $A$ is a $2p$-th root of
unity, then $$\< L\> \equiv \overline{\<L\>} \mod (p).$$  
\end{corollary}
Here and below $\bar z$ denotes the complex conjugate of $z.$ 

This last result motivates Theorem \ref{main} which gives a condition
that has to be satisfied by the quantum invariant of
a $p$-periodic pair $(M,L)$, where $M$ is a homology 3-sphere, $L$
is a colored framed link in $M$, and  $p$ is an odd prime.
On the other hand, Corollary \ref{periodicbracket2} follows from the
proof of Theorem \ref{main}
by taking $M$ to be the 3-sphere, and coloring $L$ with one.

\section{Periodic Manifolds}\label{3} 
In this section we show that
any periodic homology sphere can be presented as the result of surgery on
a periodic link. We use this to observe that any periodic homology sphere can be
obtained by doing equivariant surgery along a knot in a simple covering space of
a homology $S^1 \times S^2$.

Let $n$ be an integer greater than $1$. A $\bZ_n$ action is called
semi-free if it is free outside of the fixed point set.
\begin{definition} 
An orientable $3$-manifold $M$ is 
$n$-periodic if it admits smooth semi-free action of the cyclic group
$\bZ_n$ with a circle as the fixed point set.  
\end{definition} 

We need the following elementary lemma. First fix some notation. Suppose
that $\bZ_n$ acts on $M$ with the fixed-point set equal to a circle
$\gamma$. Denote the quotients by $M_*=M/Z_n$, and $\gamma_*=\gamma$, and
the projection map by $q:M\to M_*$. 
\begin{lemma}\label{epi} 
The map
$q_*:H_1(M)\rightarrow H_1(M_*)$ is an epimorphism. 
\end{lemma}

\begin{proof} Let $x_0\in\gamma$. Since $x_0$ is a fixed-point of the
action, any
 loop based at $q(x_0)$ lifts to a loop based at $x_0$. Thus $q_{\# }:
\pi_1(M,x_0)\rightarrow \pi_1(M_*,q(x_0))$ is an epimorphism, and since
$H_1$ is an abelianization of $\pi_1$, the map $q_*$ is also an
epimorphism. \end{proof} 

\begin{theorem}\label{periodicsurgery} 
Let $M$ be
a homology sphere. If $\bZ_n$ acts faithfully
on $M$ with the fixed-point set equal to
a circle then $M$ is $n$-periodic and
there exists an $n$-periodic framed link $L \subset S^3$ such that
$M$ is the result of surgery on $L$, and the linking number of each
component of $L$ with the fixed point set for the action implied by the
$n$-periodicity of $L$ is equal to zero. 
\end{theorem}
 
\begin{proof} 
It follows from Smith theory \cite{Br} that any
faithful
$\bZ_n$ action on a homology
sphere which has a circle as a fixed point set has to be
semi-free. Thus $M$ is $n$-periodic.

The quotient manifold  $M_*$
can be represented as a result of surgery on some framed link $L_*$ in
$S^3$. We can assume that $L_*$ satisfies the following conditions:
\begin{enumerate} 
\item $\gamma_*\cap L_*=\emptyset$; 
\item ${\rm lk}(L_*^i, \gamma_*)=0$, for any component $L_*^i$ of
  $L_*$; 
\item $\gamma_*$ is unknotted in $S^3$. 
\end{enumerate}
$L_*$ satisfying the conditions 1-3 can be obtained as follows:\\ Let
$L'_*$ be a framed link in $S^3$ such that $M_*$ is a result of surgery on
$L'_* $. Let $\hat L_*$ denote the co-core of the surgery. Co-core is a
central line, $\{0\} \times S^1$, of the solid torus, $D^2 \times S^1$,
attached to $S^3- int(V_{L'_*})$ in order to obtain $M_*$. As usual,
$V_{L'_*}$ denotes a tubular neighborhood of $L'_*$ in $S^3$. In
particular $\hat L_*$ is a framed link in $M_*$ such that $S^3$ is a
result of surgery on $\hat L_*$. Since $\gamma_*$ is homologically trivial
in $M_*$ it bounds a surface $F_*$. We can use isotopy to push $\hat L_*$
outside of $F_*$, therefore $F_*$ survives surgery. It has boundary
$\gamma_*$ and is disjoint from $L_*$ (obtained from $L'_*$) in $S^3$.
Conditions (1) and (2) follow.
 
Now we are ready to unknot $\gamma_*$ using the Kirby calculus (\cite{K},
\cite{FR}). Choose some orientation on $\gamma_*$. We can add unlinked
components with framing $\pm 1$ to $L_*$ around each crossing of
$\gamma_*$, making sure that arrows on $\gamma_*$ run opposite ways (i.e.
the linking number of $\gamma_*$ with the new component of $L_*$ is zero).
Use the K-move to change the appropriate crossings (see Figure
\ref{gamma}), and thus unknot $\gamma_*$. Thus we obtained condition (3)
without compromising conditions (1) and (2). 
\begin{figure}[ht!] \mbox{}\hfill
 \ \epsfig{file=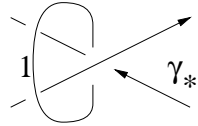, width=1in} \hfill\mbox{} 
\caption{Unknotting $\gamma$} \label{gamma}
\end{figure} 
Consider the $n$-fold cyclic branched covering of $S^3$ by
$S^3$ over $\gamma_*$. Let $L$ denote the pre-image of $L_*$. Notice that
$L$ is $n$-periodic. We claim that the result of performing surgery on $L$
is $\bZ_n$-homeomorphic to $M$. The pre-image of each component of $L_*$
consists of $n$ components permuted by $\bZ_n$ action, by condition (2).
Therefore $\bZ_n$ acts on the result of surgery on $S^3$ along $L$,
$(S^3,L)$, with branched set $\gamma_*$ and quotient $(S^3,L_*)$. By Lemma
\ref{epi} the manifold $M_*=(S^3,L_*)$ is a homology sphere. Therefore
$H_1(M_*-\gamma_*)=\bZ$, and thus $M_*-\gamma_*$ has unique $n$-fold
cyclic covering. As a consequence, the branched covering of
$(M_*,\gamma_*)$ is also unique.
Thus $(S^3,L)$ is $\bZ_n$-homeomorphic to $M$. 
\end{proof} 

A $\bZ_n$-covering space 
$X$ of  a space $X_*$
is called {\em simple} if it is classified by
a map $H_1(X_*) \rightarrow {\mathbb Z}_n$ which factors through an
epimorphism $\chi$ to ${\mathbb Z}$ \cite{Gi1,Gi2}. 
One can
define a {\em signature defect} \cite{hirz} of a simple covering $N\rightarrow N_*$ of closed orientable 3-manifolds
as follows.
Since $N\rightarrow N_*$
is a boundary of a simple $\bZ_n$ covering space of $4$-manifolds
$W\rightarrow W_*$,
we 
may
define $\text{ def }( N \rightarrow N_*)$ to be
$\text{ def }( W \rightarrow W_*)=n\sigma(W_*)-\sigma(W)$. This does
not depend on the choice of $W$.

\begin{remark} 
The third author  has
 obtained a version of Theorem 7
  for  any $p$-periodic closed oriented 3-manifold $M$, where $p$ is
  prime. This appeared 
  in \cite{PS}. The third author and Sokolov also obtained a similar
  theorem for free $\bZ_p$ actions \cite{PS}. More recently Sakuma has
  shown that every orientation preserving  finite cyclic action on a
  closed oriented 3-manifold $M$ has an equivariant framed link
  description \cite{Sa}. 
\end{remark}

\begin{corollary}\label{periodicsurgery2} 
Suppose that $M$ is an $n$-periodic homology sphere. Then $M$ can be
obtained by equivariant framed surgery along a knot in a simple covering of a homology
$S^1\times S^2$. Moreover, the signature defect of such a covering is
always even.
\end{corollary}

\begin{proof}
Suppose that $\bZ_n$ acts on $M$ with the fixed-point set equal to a
circle $\ga$.
Let $L$ be a periodic framed link surgery description of $M$ given by
Theorem \ref{periodicsurgery}, and let $N$ be given by surgery on the
framed link $L_+= L\cup \gamma$, where $\gamma$ has framing zero
(i.e.,  $N$ is obtained by performing zero-framed surgery on $M$ along
$\ga$). 

The free action on $M\setminus \ga$ extends to a free $\bZ_n$-action on $N$.
Denote by $N_*$  the orbit space of this action. The manifold $N_*$ can be obtained by
doing zero framed surgery on $M_*$ along $\ga_*$. (We keep the
notation from the proof of Theorem \ref{periodicsurgery}.)
 By Lemma \ref{epi}, $M_*$
is a homology sphere, so $N_*$ is a homology $S^1 \times S^2.$ 
It follows 
that the quotient map $N \ra N_*$ is a simple $\bZ_n$-covering space. 
By \cite[Prop 4]{Gi1}, $\text{ def }( N \rightarrow N_*)$ is even.

Let $\delta$ be a meridian of $\ga$ in $M$. The curve $\delta$ bounds a
disk transverse to $\ga$. We frame $\delta$ so its push-off
bounds a transverse disk to $\gamma$ which misses the original
transverse disk. If $\delta'$ denotes the image of $\delta$ in $N$
under the surgery on $M$ along $\ga$ then $M$ can be recovered from
$N$ by doing surgery along $\delta'$.
\end{proof}

If $L$ is a framed link in $S^3,$ the linking matrix of
$L$ has off-diagonal entries given by the pairwise linking numbers of the
components, and diagonal entries given by the linking numbers of the
components with their push-offs according to the framing. 
Recall that $\sigma (L)$  denotes the signature of the linking matrix of $L$.
\begin{remark}\label{signatures}
With the notation from the proof above, note that 
$L_{++}= L \cup \ga \cup \de$ is another framed link surgery description of $M$.
 We have
that $\sigma(L)=\sigma(L_+)= \sigma(L_{++})$.
\end{remark}

\section{ $SO(3)$-Quantum Invariants at odd primes.}\label{4} 

Throughout the rest of this paper we assume $p$ is an odd prime.  
Let $M$ be a closed, connected, oriented 3-manifold, and $J$ a framed
oriented $p$-colored  
link in $M$, i.e., a link whose components have been
assigned nonnegative integers less than $p-1$. 
We begin this section by recalling the definition of the $SO(3)$ invariant
${\mathcal I}_p(M,J)$, following  the description in \cite{MR}.

Suppose that $M$ is represented by a
surgery on a framed link $L\subset S^3,$ and that $J$ is disjoint from the
co-cores of the surgery, so  $J$ lies in $S^3 \setminus L.$ We
use the renormalized Kauffman bracket 
\[ [L]=-\left(A^2+A^{-2}\right)\langle L\rangle .\] 
Notice that $[\emptyset]=1$. Set
the variable $A$ to be a primitive root of unity of order $2p$. We denote
by $k$ a root of unity such that $$k^2 = A^{-6 -p(p+1)/2}.$$ 
If $p\equiv -1\pmod{4}$ then $\bZ[A,k]= \bZ[A].$ If $p\equiv
1\pmod{4}$ then $\bZ[A,k]= \bZ[A,i]$. 
Elements $e_i$ of the skein module of a solid torus are defined
recursively by $e_0=1$, $e_1=z$ and 
\[e_{i+1}=ze_i-e_{i-1},\] 
where $z$ is represented by a longitude of the torus.  
The value of the Kauffman bracket
of the skein element $e_i$, when the solid torus is embedded
in $S^3$ in a standard way, is given by 
$$\Delta_i = (-1)^i\frac{A^{2i+2}-A^{-2i-2}}{A^2-A^{-2}}.$$
By $\Omega_p$ 
we denote the following element of the Kauffman bracket skein
module of a solid torus:
$$\Omega_p=\sum_{i=0}^{\frac {p-3} 2} \Delta_{i} e_{i}.$$ 
Let 
$$\eta = k \frac 1{2p} A^3 (A^2-A^{-2} )A^{\frac {p (p-1)}2}\sum_{m=1}^{2p}(-1)^m A^{m^2}.$$
which satisfies
$$\eta^2[\Omega_p] =1, \quad \text{i.e} \quad \eta ^2= \frac{-(A ^2-A
^{-2})^2}{p}.$$ 
Recall that the invariant ${\mathcal I}_p(M,J)$ is given by:
\begin{equation}\label{niezmiennik} 
{\mathcal I}_p(M,J)=k^{-\sigma_L}[
L(\eta \Omega_p) \sqcup J]. 
\end{equation} 
Here $[ L(\eta \Omega_p) \sqcup J]$
 denotes the bracket of the skein obtained by replacing each component of
$L$ by $\eta \Omega_p,$ and replacing each component of $J$ by the skein
$e_i,$ where $i$ is the color of the component.
If one sets 
$A=e^{2\pi i (1+p^2)/4p}$, and chooses $k$ so that  $\eta >0$, the invariant ${\mathcal I}_p$ is
equal to the quantum $SO(3)$-invariant $\tau_p^{SO(3)}$ \cite{KM}. 

Our
criteria for periodicity of links require that we consider invariants
modulo $p$. 
 From the definition it is only apparent  that
${\mathcal I}_p(M,J)\in \bZ[A, k , \frac 1 p].$ Thus we will need
the following result proved by Murakami \cite{Mu} when $J$ is empty
and generalized to colored links by Masbaum and Roberts \cite{MR}.  
\begin{theorem}\label{Murakami} 
Let $J$ be a framed
$p$-colored link in  $M,$ 
then ${\mathcal I}_p(M,J)\in \bZ[A,k]$. 
If $p\equiv -1 \pmod 4$ or $ \beta_1(M)$ is even (in particular if $M$ is a
homology sphere) then 
${\mathcal I}_p(M,J)\in \bZ[A].$ 
\end{theorem}

 In the proof of our criterion we also use the following integrality result for
 closed $3$-manifolds with non-zero first Betti number. It was proved
 by Murakami \cite{Mu}, and by Cochran and Melvin \cite{CM} for empty
 $J$. In full generality it follows from Gilmer's \cite[Corollary 4]{Gi1}.

We adopt the notation: 
$$\<\<M,J\>\>_p=\eta {\mathcal I}_p(M,J).$$ 
This is suggested by 
 TQFT notation (see the Appendix).

\begin{theorem}\label{betti} 
Suppose that $M$ is connected and has first Betti
number non\-zero, and $J$ is a framed $p$-colored link in $M,$ then
$\<\<M,J\>\>_p \in \bZ[A,k]$. 
\end{theorem}

\section{$SO(3)$-Quantum Invariants and $\bZ_p$-group actions}\label{5} 

We will show in this
section how the $p$-periodicity of $M$ is reflected in quantum
$SO(3)$-invariant at $p,$ assuming that  $p$ is an odd prime. 
We will need the following theorem of Gilmer \cite{Gi1}.   
For  the convenience of the reader, a proof of this theorem
is sketched in an appendix.
Notice, that if $N$ is a simple $\bZ_{p}$ covering space of $N_*,$ then
both $N$ and $N_*$ have first Betti number nonzero.  Thus, by
Theorem \ref{betti}, it makes sense to consider $\<\<N,J\>\>_p$ modulo $p$.

\begin{theorem} \label{simplecong} 
If $(N,J)$ is a simple $\bZ_{p}$ covering
space of $(N_*,J_*)$,  where $N_*$ is connected with odd first Betti number, and
 $p$ is an odd prime, then for some $n \in \bZ,$ we have: $$\<\<N,J\>\>_p = k
^{\text{ def }( N \rightarrow N_*)} n \mod (p).$$ 
\end{theorem}

We are ready to state and prove our main theorem. We remark that this theorem is uninteresting in the case $p=3$ as the conclusion always holds whether $M$ is 3-periodic or not.

\begin{theorem}\label{main} 
Suppose that $M$ is a homology sphere, and $J$
is a  $p$-colored framed link in $M.$ If $(M,J)$ is $p$-periodic and
$p$ is an odd prime then 
$${\mathcal I}_p(M,J) = A^{2j}\cdot \overline{{\mathcal I}_p(M,J) }\mod (p)$$
 for some integer $j$. 
\end{theorem} 
\proof 
We are using the notation introduced in the proof of Corollary
\ref{periodicsurgery2}.
By Remark
\ref{signatures}, $M$ is the result of surgery along $L_{++}$, so:
$${\mathcal I}_p(M,J) = k^{-\sigma(L_{++})}[ L_{++}(\eta \Omega_p) \sqcup J].$$
Expanding $\delta$ cabled with $\eta \Omega_p,$ where $\delta_i$ denotes
$\delta$ colored with $i,$ we obtain: 
$${\mathcal I}_p(M,J)= k^{-\sigma(L_{+})} \eta
\sum_{i=0}^{\frac {p-3} 2} \Delta_i [ L_{+}(\eta \Omega_p) \sqcup J \sqcup \delta_i],$$
 which is equal to
$$\sum_{i=0}^{\frac {p-3} 2} \Delta_i \left(\eta k^{-\sigma(L_{+})} [L_{+}(\eta \Omega_p) \sqcup J
\sqcup \delta_i]\right)$$ 
or, by the definition of $\<\<\quad \>\>,$ 
$$ \sum_{i=0}^{\frac {p-3} 2} \Delta_i \<\<N,J\sqcup \delta_i\>\>.$$
Theorem \ref{simplecong}  implies that, modulo $p$, this is
$$ \sum_{i=0}^{\frac {p-3} 2} \Delta_i k^{\text{ def }( N \rightarrow N_*)}n_i,$$ 
where $n_i \in \bZ.$ 
Since, by Corollary \ref{periodicsurgery2}, $\text{ def }( N \rightarrow N_*)$ is even,
$${\mathcal I}_p(M,J)= k^{2s} \sum_{i=0}^{\frac {p-3} 2} n_i \Delta_i \pmod p,$$ 
where $s \in \bZ.$ 
Noting that $\Delta_i=\overline{ \Delta_i}$, we have that 
$${\mathcal I}_p(M,J)=k^{4s} \overline{ {\mathcal I}_p(M,J)} \pmod p.\eqno{\qed}$$

The proof shows that phase
factor is $k^{2 \text{ def }( N \rightarrow N_*)}.$ If $M$ is $S^3,$ this
defect is zero, and we obtain Corollary \ref{periodicbracket2}, as a
consequence of the proof of Theorem \ref{main}. 

One expects similar results for quantum invariants
associated to other Lie groups. These results should be related to
Traczyk-Yokota periodicity criteria involving skein polynomials \cite{P2},
\cite{T2}, \cite{Y}. Corollary 5 in \cite{Gi1} is a first step in this
direction.  
In fact, subsequent to an earlier version of this paper, the results of this paper have been generalized
to quantum invariants associated to semi-simple Lie groups by Chen and Le \cite{CL}.

There should be a criterion for periodicity of other $3$-manifolds
than homology spheres in terms of their $SO(3)$-invariants.
At present,  we do not know how to prove it since
 no version of Theorem  \ref{simplecong} has  been obtained
for non-simple $\bZ_p$ covers.   

\section{Examples}\label{6} 
\begin{figure}[ht!] 
\mbox{}\hfill \
\epsfig{file=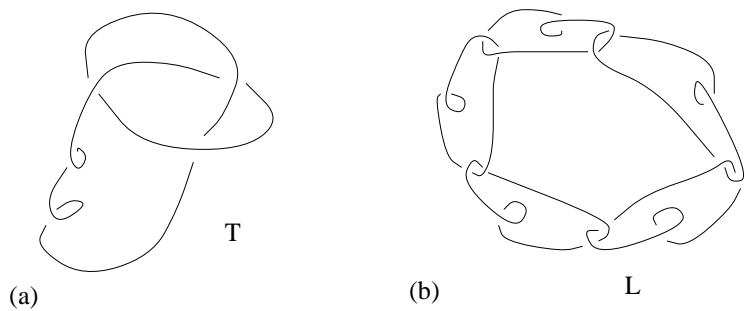,height=4cm} \hfill\mbox{} 
\caption{Surgery descriptions of the Poincar\'e homology sphere}
\label{PHS} 
\end{figure}

We illustrate Theorem \ref{main} by
considering the Poincar\'e homology sphere $\Sigma$ \cite{Po}. It is well
known that $\Sigma$ is $5$-periodic since it is the 5-fold branched cyclic
cover of $S^3$ along the trefoil knot. The $5$-periodicity is also
apparent by presenting $\Sigma$ as a result of surgery on a $5$-periodic
link $L$ pictured in Figure \ref{PHS} (compare \cite{Gol}, \cite{Ro}) with
blackboard framing. 
One can compute ${\mathcal I}_5(\Sigma)$ from its description as $+1$
surgery on the right handed trefoil \cite{D} (knot $T$ in Figure
\ref{PHS}):
$${\mathcal I}_5(\Sigma)= 1- 2A + 2A^2 -A^3.$$ 
Thus 
$${\mathcal I}_5(\Sigma)= A^{-2}(1 - 2A^{-1} + 2A^{-2} - A^{-3})
= A^{-2}\overline{{\mathcal I}_5(\Sigma)}$$ 
i.e. $j=-1$ in
our criterion. It also follows from Theorem \ref{main} that $\Sigma$
is not $7$-periodic. We computed ${\mathcal I}_7(\Sigma)$ from the
same surgery presentation. 
$$
{\mathcal I}_7(\Sigma)= -2+A+2A^3-A^4 
$$
One can check that 
$$
{\mathcal I}_7(\Sigma)\not\equiv A^{2j} \overline{{\mathcal I}_7(\Sigma)} \pmod{7}
$$
for any $j$. Thus $\Sigma$ is not $7$-periodic.

More generally, consider the Brieskorn homology spheres 
$\Sigma(2, |n|, |2n-1|),$ for $n$ odd.  Let $K_n$ be a right handed $(2,n)$
torus knot if $n$ is positive, and a left handed $(2,-n)$ torus knot if $n$
is negative. Denote by $M_ n$  the manifold obtained by $+1$ surgery
along $K_n.$  It was shown by Seifert in \cite{S} that
$M_n$ is a Seifert fibered homology sphere, with 
three exceptional fibers with multiplicities $2,$  $|n|,$ and
$|2n-1|.$  Moreover, also according to  \cite{S}, fixing $n$, there is
only one  homology sphere with
such Seifert fibers (up to orientation). Thus $M_ n$ is, up to
orientation, the Brieskorn homology sphere $\Sigma(2, |n|,|2n-1|)$.
One can evaluate ${\mathcal I}_p(M_n)$ using the recoupling
theory of \cite[Chapter 9]{KL}.
To obtain the value of the bracket of $K_n$ colored with $i$, use the
first part of (9.15)  to fuse the two strands, then apply (9.9) to
untwist the resulting theta curve, and finally  evaluate the resulting
planar  theta curve using (9.10). We obtain: 
$${\mathcal I}_p(M_n)=\eta k^{-1}\sum_{i=0}^{\frac {p-3} 2} 
\left( (-A)^{i (i+2)(1-n)} \Delta_i\sum_{j=0}^i
  \Delta_{2j}(\lambda_{2j}^{i \  i})^{-n}\right), $$ 
where \cite[9.9]{KL}:
$$\lambda_{2j}^{i \  i}= (-1)^{i-j} A^{i(i+2)-j(2j+2)}.$$
Lawrence derived the SU(2) invariants of this manifold in a  related
way \cite{L}. 

We programmed Mathematica to evaluate ${\mathcal I}_p(M_n)$ symbolically
and check\hfill \linebreak
whether the conclusion of Theorem \ref{main} holds.
Since  $M_n$ is  $p$-periodic for any odd prime $p$  dividing 
$n$ or $2n-1,$  we know the conclusion should hold in this case. 
We checked to see if the conclusion held for odd primes $5\le p\le 19, $
and for
 odd integers satisfying
$3 \le |n| \le 19.$ The condition held as expected whenever $p$
divided 
$n$ or $2n-1.$ This happened in 
$27$
cases. It also held for the following
list of  
ten
other pairs of $(n,p):$ 
$(9,5)$, $(11,5)$, $(13,7)$, $(15,7)$, 
 $(19,5)$,  
$(-9,5)$, $(-11,5)$, $(-13,7)$, $(-15,7)$, $(-19,7)$. 
For the 
$71$
other pairs in the above  range, we can say $\Sigma(2, |n|, |2n -
1|)$ is not $p$-periodic. Thus we checked the conclusion for  
$6\times
2\times 
9=108$
different
pairs of prime $p$ and manifold $M$. The conclusion of Theorem
\ref{main} held in only
$37$ cases. In 
$27$ of these cases the manifold is known to be
$p$--periodic,
as mentioned above.

In addition, we checked the criterion  for the Poincar\'e homology
sphere $\Sigma=M_3 $ for larger values of $p$, and concluded that 
$\Sigma$ is not $p$-periodic if $p$ is any  prime between $7$ and $61$
(where we stopped the computer). 

\begin{remark}
Zhang has proven that the only periods for the Poincar\'e homology sphere are
$2$, $3$, and $5$, \cite{Zh}.
More generally, we were informed by  Boileau \cite{Bo}
that it follows from the orbifold theorem \cite{Bo-Po} and the previous work
of Meeks and Scott \cite{M-S} and Seifert that if 
the Brieskorn homology sphere, $\Sigma (2, |n|, |2n-1|)$ ($|n|>2$), is
$Z_p$ periodic ($p$ odd prime) then $p$ divides $n$ or $2n-1$.
\end{remark}

\begin{remark}
Chen and Le \cite{CL} checked that our periodicity criterion is also sufficient
to prove that the only possible periods for the Poincar\'e homology
sphere are $2$, $3$ and $5$, and the only possible periods for
$\Sigma(2,3,7)$ are $2$, $3$ and $7$.

\end{remark}

\Addresses\recd\newpage

\cl{\large\bf Appendix: Quantum invariants of simple free $\bZ_p$ covers}
\bigskip

\small
\cl{\sc Patrick M. Gilmer}\bigskip

As background for the proof of Theorem \ref{simplecong}  we now
describe some features of the topological quantum field theory 
defined in \cite{BHMV2}, and the
relation of the functor $(V_p,Z_p)$ to  
${\mathcal I}_p(M,J),$ and $\<\<M,J\>\>_p.$

The notion of a $p_1$-structure  provides one solution to the problem
of finding extra structures for surfaces and 3-manifolds so that one
can remove the ``framing anomalie'' in the Reshetikhin-Turaev-Witten
TQFTs.  
See \cite{BHMV2} for the definition of  $p_1$-structures.
 
For a closed 3-manifold $M$ with $p_1$-structure $\alpha,$ there is an
associated integer $\sigma(\alpha)$ which specifies the
$p_1$-structure up to homotopy. 

Let $\kappa$ be a cube root of $k.$ Let $R_p$ be the ring
$\bZ[A,\kappa, \frac 1 p].$ 
The functor $V_p$ assigns a free finitely generated module over $R_p$ to a
surface $\Sigma$ equipped with a $p_1$-structure and some even-colored
framed points. 
This module is denoted $V_p(\Sigma).$
Let $N$ be a 3-dimensional cobordism with a framed even-colored link J
and $p_1$-structure,  whose boundary is identified with $-\Sigma_1 \sqcup
\Sigma_2$.  Link $J$ must meet the boundary in the colored framed points. 
 The functor $Z_p$ assigns to $N$ a 
module homomorphism  $Z_p(N):V_p(\Sigma_1) \rightarrow V_p( \Sigma_2).$

A closed 3-manifold $M$ with $p_1$-structure and an even-colored framed link $J$ can be considered
as a cobordism from the empty set to the empty set. As $V_p(\emptyset)= R_p,$ 
$Z_p(M)$ is multiplication by some scalar in $R_p.$ This scalar is denoted
$\<M\>_p$.   The invariant ${\mathcal I}_p(M,J)$ is a rescaled version of the
invariant $\<M\>_p$ with the dependence on the $p_1$-structure
removed. One has: 
$$\<M\>_p = \kappa^{\sigma({\alpha})}   \eta {\mathcal I}_p(M,J),$$
and
$$\<\<M,J\>\>_p =\kappa^{-\sigma({\alpha})}\<M\>_p.$$  
Thus $\<\<M,J\>\>_p$ is $\<\tilde M,J\>_p,$
where $\tilde M$ is $M$ with a reassigned $p_1$-structure whose
$\sigma$ invariant is zero. 
Note that when we use
the quantum invariant notations: ${\mathcal I}_p$ and $\<\<\quad\>\>_p, $
we explicitly include $J.$ 
However when using the TQFT notations: $Z_p$ and $\< \quad\>_p,$ the
symbols for the cobordisms $M$ and $N$ etc., implicitly include the
data of $J$ and a $p_1$ structure. 

If $\sigma(\alpha)= 3(\beta_0(M)+ \beta_1(M)) \pmod 6,$  one has 
\cite[\S 2]{BHMV2}\cite{BHMV1} that
$\<M,J\>_p  \in \bZ[A, \frac 1 p].$

Suppose that $E$ is a cobordism from  $\Sigma$ to itself. Then one may
form $C(E)$, the closure of  $E$, by identifying the two copies of
$\Sigma$ which make up the boundary.  
We have the important trace property:
$$\<C(E)\>_p= \text{Trace} (Z_p(E)).$$
If $H$ is a connected cobordism from $\emptyset$ to $\Sigma,$ let [H]
denote the image of  
$1 \in V_p(\emptyset)$ in $V_p(\Sigma)$ under $Z_p(H).$ The module
$V_p(\Sigma)$ is finitely generated by such elements. Moreover there
is a well-defined non-degenerate  
Hermitian form $\<\quad ,\ \>_{\Sigma} $ on  $V_p(\Sigma)$ given by 
$$\<[H_1],[H_2]\>_{\Sigma}=\<H_1 \cup_{\Sigma} -H_2\>_p.$$
This completes the summary of the properties of $(V_p,Z_p)$ that we need
for the proof of  Theorem \ref{simplecong}.

\proof

If $N_*$ and $N$ contain some link components with odd colors, trade
them for even colors using the fact that exchanging
 $e_{\frac {p-1}{2}+i}$ with $ e_{\frac {p-3}{2}-i}$  
 does not  change the value of $\<\<N\>\>_p$ \cite{BHMV1}.
Give $N_*$ a $p_1$-structure $\alpha_*,$ with $\sigma(\alpha_*)= 0,$
so that  $\<N_*\>_p  \in \bZ[A, \frac 1 p]$.  By  Theorem \ref{betti}
$\<\<N_*,J_*\>\> _p \in \bZ[A,k].$ Thus 
$\<N_*\>_p \in \bZ[A,k].$ Therefore $\<N_*\>_p \in \bZ[A,k] \cap
\bZ[A, \frac 1 p]=\bZ [A].$ 

$N$ acquires an induced $p_1$-structure $\alpha$. 
By \cite[Lemma 4]{Gi2}, we have that 
$3 \text{ def }( N \rightarrow N_*)= -\sigma(\alpha).$ 
Since  $k=\kappa^3,$ we have that 
$\<\<N,J\>\> =k^{ \text{ def }( N  \rightarrow N_*)} \<N\>_p.$   

Since $N \rightarrow N_*$ is simple, it can be defined using reduction
modulo $p$  of some epimorphism $\chi:H_1(N_*) \rightarrow Z$.
Let $\Sigma$ be a connected surface in $N_*$ which is Poincare dual to
the class of $\chi$ in $H^1(N_*).$  
Let $E$ denote $N_*$ slit along $\Sigma.$ 
The
manifold $E$ is
connected and we may regard $E$ as a cobordism from $\Sigma$ to $\Sigma.$
We denote the composition in the cobordism category  of $E$ with itself $n$ times by $E^n.$ 
We have that $N_*= C(E),$ and $N$ is diffeomorphic to $C(E^p).$ 
Thus $\<N_*\>_p= \text{Trace}(Z_p(E)),$ and $\<N\>_p= \text{Trace}((Z_p(E))^p),$

We will show below that there is a matrix ${ \mathcal {M}}$ whose
coefficients are algebraic integers such that 
$ \text{Trace}({  \mathcal {M}}) =\text{Trace}(Z_p(E)),$ 
and $\text{Trace}({ \mathcal {M}}^p)=\text{Trace}(Z_p(E)^p).$ 
Assuming the  existence of ${ \mathcal {M}},$  we now complete the proof. 

We  have that  $\text{Trace}({\mathcal {M}})= \sum_i \lambda_i,$   where
$\{\lambda_i\}$ denotes the list
of
eigenvalues of ${\mathcal {M}},$ with
multiplicities. Each $\lambda_i$ is  an algebraic integer as it
satisfies  the characteristic
polynomial of ${\mathcal{M}},$  whose coefficients are algebraic integers. Moreover $\text{Trace}({\mathcal {M}}^p)=
\sum_i \lambda_i^p.$ Moreover   $\sum_i \lambda_i^p-(\sum_i
\lambda_i)^p$ is $p$ times an algebraic integer, since $\binom{p}{j} =0 \pmod{p} $  for $0<j<p.$

So $ \<N\>_p -(\<N_*\>_p)^p$ is $p$ times an algebraic integer. As $\<N_*\>_p \in \bZ [A],$
we have that $(\<N_*\>_p)^p$ is congruent modulo $p$ to an ordinary
integer, say $n.$ Thus 
 $\<N\>_p -n $ is $p$ times an algebraic integer. It follows that 
$\<\<N,J\>\> -k^{\text{ def }( N \rightarrow N_*)}  n$  is $p$ times an algebraic integer.
This is the conclusion of the theorem.

Thus we only need to find ${ \mathcal {M}}.$ 
Let $S_p(\Sigma)$ denote the $\bZ[A,\kappa]$ submodule of $V_p(\Sigma)$ generated
by all elements of the form $[H]$ where $H$ is a connected cobordism from $\emptyset$ to $\Sigma.$ 
As this set spans $V_p(\Sigma)$ over $R_p,$ we have that 
$S_p(\Sigma) \otimes_{\bZ[A,\kappa]} R_p= V_p(\Sigma).$

We claim that $S_p(\Sigma)$ is  a finitely generated projective 
$\bZ[A,\kappa]$ -module. To see this
  pick a finite set of generators:
  $\{ [H_{i}] \}_{1 \le i\le m}$ for $V_p(\Sigma)$  where
 $H_{i}: \emptyset  \rightarrow \Sigma$ are connected.
 We have an injective map
$$S_p(\Sigma) \rightarrow \bZ[A,\kappa]^{m},$$ given by
$$[H] \mapsto \ {\frac 1 \eta}  ( \<[H],[H_{1}]\>_{\Sigma}, \cdots,
\<[H],[H_{m}]\>_{\Sigma} ).$$ 
The reason that this map takes values in 
$\bZ[A,\kappa]^{m}$ is Theorem \ref{Murakami}.
As $\bZ[A,\kappa]$ is a Dedekind domain, a submodule of
$\bZ[A,\kappa]^{m}$ is necessarily finitely generated and projective
\cite[Proposition 10.12]{Ja}.

The $\bZ[A,\kappa]$ submodule $S_p(\Sigma)$ of $V_p(\Sigma)$ is preserved by
$Z_p(E):V_p(\Sigma) \rightarrow V_p(\Sigma).$ This is seen as follows.
Given  a connected $H: \emptyset \rightarrow \Sigma,$ $Z_p(E)[H]=
[H \cup_{\Sigma} E],$ and $H \cup_{\Sigma} E$ is connected. We denote the induced map by 
${\mathcal{Z}}:S_p(\Sigma) \rightarrow S_p(\Sigma).$

There is an ideal $\mathcal{J}$ in $\bZ[A,\kappa]$ such that 
$S_p(\Sigma) \oplus \mathcal{J} \approx \bZ[A,\kappa]^l$, where
$l=\text{dim}V_p(\Sigma)+1$ 
\cite[Theorem 10.14]{Ja}.
So, denoting by $ \text{0}_{\mathcal{J}}$ the zero map from
$\mathcal{J}$ to itself,
$$\text{Trace}(Z_p(E))= \text{Trace}( {\mathcal{Z}}\oplus
\text{0}_{\mathcal{J}})= 
\text{Trace}( \mathcal {M}),$$  
where ${ \mathcal {M}}$ is a matrix which represents the endomorphism 
${\mathcal {Z}}  \oplus 0_{\mathcal{J}}$ of free $\bZ[A,\kappa]$
modules. Similarly
$$\text{Trace}(Z_p(E^p))= \text{Trace}( {\mathcal {Z}}^p \oplus 0_{\mathcal{J}}) =
\text{Trace}( {\mathcal {M}}^p).\eqno{\qed}$$

\end{document}